\documentstyle[amsmath,amssymb]{article}
\newtheorem{th}{Theorem}
\newtheorem{lem}{Lemma}

\newtheorem{cor}{Corollary}
\newtheorem{defn}{Definition}
\newenvironment{defn-new}{\begin{defn} \em}{\end{defn}}
\newtheorem{rem}[th]{Remark}
\newenvironment{rem-new}{\begin{rem} \em}{\end{rem}}
\newtheorem{ex}[th]{Example}
\newenvironment{ex-new}{\begin{ex} \em}{\end{ex}}
\newtheorem{exer}[th]{Exercise}
\newenvironment{exer-new}{\begin{exer} \em}{\end{exer}}
\newtheorem{agr}[th]{Agreement}
\newenvironment{agr-new}{\begin{agr} \em}{\end{agr}}
\newtheorem{pbm}[th]{Problem}
\newenvironment{pbm-new}{\begin{pbm} \em}{\end{pbm}}

\makeatletter \@addtoreset{equation}{section} \makeatother

\begin{document}

\begin{center}
{\Large {\bf Index of quasi-conformally symmetric semi-Riemannian manifolds}}
\end{center}
\begin{center}
{\large Mukut Mani Tripathi, Punam Gupta, Jeong-Sik Kim}
\end{center}

\medskip

\noindent {\bf Abstract.}
We find the index of $\widetilde{\nabla}$-quasi-conformally
symmetric and $\widetilde{\nabla}$-concircularly symmetric
semi-Riemannian manifolds, where $\widetilde{\nabla}$ is metric
connection.
\medskip

\section{Introduction\label{ric-sym-hay-con-index}}

In 1923, L.P. Eisenhart \cite{Eisenhart-23} gave the condition for
the existence of a second order parallel symmetric tensor in a
Riemannian manifold. In 1925, H. Levy \cite{Levy-25}\ proved that a
second order parallel symmetric non-singular tensor in a real space
form is always proportional to the Riemannian metric. As an
improvement of the result of Levy, R. Sharma \cite{Sharma-89}\
proved that any second order parallel tensor (not necessarily
symmetric) in a real space form of dimension greater than $2$ is
proportional to the Riemannian metric. In 1939, T.Y. Thomas \cite
{Thomas-39} defined and studied the index of a Riemannian manifold.
More precisely, the number of metric tensors (a metric tensor on a
differentiable manifold is a symmetric non-degenerate parallel
$(0,2)$ tensor field on the differentiable manifold) in a complete
set of metric tensors of a Riemannian manifold is called the index
of the Riemannian manifold \cite[p. 413]{Thomas-39}. Thus the
problem of existence of a second order parallel symmetric tensor is
closely related with the index of Riemannian manifolds. Later, in
1968, J. Levine and G.H. Katzin \cite{Levine-Katzin-68} studied the
index of conformally flat Riemannian manifolds. They proved that the
index of an $n$-dimensional conformally flat manifold is $n(n+1)/2$
or $1$ according as it is a flat manifold or a manifold of non-zero
constant curvature. In 1981, P. Stavre \cite{Stavre-81} proved that
if the index of an $n$-dimensional conformally symmetric Riemannian
manifold (except the four cases of being conformally flat, of
constant curvature, an Einstein manifold or with covariant constant
Einstein tensor) is greater than one, then it must be between $2$
and $n+1$. In 1982, P. Starve and D. Smaranda
\cite{Stavre-Smaranda-82} found the index of a conformally symmetric
Riemannian manifolds with respect to a semi-symmetric metric
connection of K. Yano \cite{Yano-70}. More precisely, they proved
the following result: \textquotedblleft Let a Riemannian manifold be
conformally symmetric with respect to a semi-symmetric metric
connection $\overline{\nabla }$. Then (a) the index
$i_{\overline{\nabla }}$ is $1$ if there is a vector field $U$
such that $\overline{\nabla }_{U}\overline{E}=0$ and $\overline{\nabla }_{U}%
\overline{r}\neq 0$, where $\overline{E}$ and $\overline{r}$ are the
Einstein tensor field and the scalar curvature with respect to the
connection $\overline{\nabla }$, respectively; and (b) the index $i_{%
\overline{\nabla }}$ satisfies $1<i_{\overline{\nabla }}\leq n+1$ if $%
\overline{\nabla }\,\overline{E}\neq 0$".

A real space form is always conformally flat and a conformally flat manifold
is always conformally symmetric. But the converse is not true in both the
cases. On the other hand, the quasi-conformal curvature tensor \cite%
{Yano-Sawaki-68} is a generalization of the Weyl conformal curvature tensor
and the concircular curvature tensor. The Levi-Civita connection and
semi-symmetric metric connection are the particular cases of the metric
connection. Also, a metric connection is Levi-Civita connection when its
torsion is zero and it becomes the Hayden connection\cite{Hayden-32} when it
has non-zero torsion. Thus, metric connections include both the Levi-Civita
connections and the Hayden connections (in particular, semi-symmetric metric
connections).

Motivated by these circumstances, it becomes necessary to study the index of
quasi-conformally symmetric semi-Riemannian manifolds with respect to any
metric connection. The paper is organized as follows. In Section~\ref%
{sec-index}, we give the definition of the index of a semi-Riemannian
manifold and give the definition and some examples of the Ricci-symmetric
metric connections $\widetilde{\nabla }$. In Section \ref{sec-qcc-flat}, we
give the definition of the quasi-conformal curvature tensor with respect to
a metric connection $\widetilde{\nabla }$. We also obtain a complete
classification of $\widetilde{\nabla }$-quasi-conformally flat (and in
particular, quasi-conformally flat) manifolds. In Section \ref{sec-qcc-index}%
, we find out the index of $\widetilde{\nabla }$-quasi-conformally symmetric
manifolds and $\widetilde{\nabla }$-concircularly symmetric manifolds. In
the last section, we discuss some of applications in theory of relativity.

\section{Index of a semi-Riemannian manifold\label{sec-index}}

Let $M$ be an $n$-dimensional differentiable manifold. Let $\widetilde{%
\nabla }$ be a linear connection in $M$. Then torsion tensor $\widetilde{T}$
and curvature tensor $\widetilde{R}$ of $\widetilde{\nabla }$ are given by
\[
\widetilde{T}\left( X,Y\right) =\widetilde{\nabla }_{X}Y-\widetilde{\nabla }%
_{Y}X-[X,Y],
\]%
\[
\widetilde{R}(X,Y)Z=\widetilde{\nabla }_{X}\widetilde{\nabla }_{Y}Z-%
\widetilde{\nabla }_{Y}\widetilde{\nabla }_{X}Z-\widetilde{\nabla }_{[X,Y]}Z
\]%
for all $X,Y,Z\in {\mathfrak X}(M)$, where ${\mathfrak X}(M)$ is the
Lie algebra of vector fields in $M$. By a semi-Riemannian metric
\cite{ONeill-83} on $M$,
we understand a non-degenerate symmetric $\left( 0,2\right) $ tensor field $%
g $. In \cite{Thomas-39}, a semi-Riemannian metric is called simply a metric
tensor. A positive definite symmetric $\left( 0,2\right) $ tensor field is
well known as a Riemannian metric, which, in \cite{Thomas-39}, is called a
fundamental metric tensor. A symmetric $\left( 0,2\right) $ tensor field $g$
of rank less than $n$ is called a degenerate metric tensor \cite{Thomas-39}.
\medskip

Let $\left( M,g\right) $ be an $n$-dimensional semi-Riemannian manifold. A
linear connection $\widetilde{\nabla }$ in $M$ is called a metric connection
with respect to the semi-Riemannian metric $g$ if $\widetilde{\nabla }g=0$.
If the torsion tensor of the metric connection $\widetilde{\nabla }$ is
zero, then it becomes Levi-Civita connection $\nabla $, which is unique by
the fundamental theorem of Riemannian geometry. If the torsion tensor of the
metric connection $\widetilde{\nabla }$ is not zero, then it is called a
Hayden connection \cite{Hayden-32,Yano-82}. Semi-symmetric metric
connections \cite{Yano-70} and quarter symmetric metric connections \cite%
{Golab-75} are some well known examples of Hayden connections.

Let $\left( M,g\right) $ be an $n$-dimensional semi-Riemannian manifold. For
a metric connection $\widetilde{\nabla }$ in $M$, the curvature tensor $%
\widetilde{R}$ with respect to the $\widetilde{\nabla }$ satisfies the
following condition%
\begin{equation}
\widetilde{R}\left( X,Y,Z,V\right) +\widetilde{R}\left( Y,X,Z,V\right) =0,
\label{eq-curv-con-semi-hay-1}
\end{equation}%
\begin{equation}
\widetilde{R}\left( X,Y,Z,V\right) +\widetilde{R}\left( X,Y,V,Z\right) =0
\label{eq-curv-con-semi-hay-2}
\end{equation}%
for all $X,Y,Z,V\in {\mathfrak X}(M)$, where
\begin{equation}
\widetilde{R}\left( X,Y,Z,V\right) =g(\widetilde{R}\left( X,Y\right) Z,V).
\label{eq-curv-con-hay-sym-6}
\end{equation}%
The Ricci tensor $\widetilde{S}$ and the scalar curvature $\widetilde{r}$\
of the semi-Riemannian manifold with respect to the metric connection $%
\widetilde{\nabla }$ is defined by
\[
\widetilde{S}\left( X,Y\right)
=\sum_{i=1}^{n}\varepsilon_{i}\widetilde{R}\left(
e_{i},X,Y,e_{i}\right) ,
\]%
\[
\widetilde{r}=\sum_{i=1}^{n}\varepsilon_{i}\widetilde{S}\left(
e_{i},e_{i}\right) ,
\]%
where $\{e_{1},\ldots ,e_{n}\}$ is any orthonormal basis of vector
fields in the manifold $M$ and $\varepsilon_{i}=g(e_{i},e_{i})$. The
Ricci operator $\widetilde{Q}$\ with respect to the metric
connection $\widetilde{\nabla }$ is defined by
\[
\widetilde{S}\left( X,Y\right) =g(\widetilde{Q}X,Y),\qquad X,Y\in {\mathfrak X}%
(M).
\]%
Define
\begin{equation}
\widetilde{e}\,X=\widetilde{Q}X-\frac{\widetilde{r}}{n}X,\qquad X\in {\mathfrak X}%
(M),  \label{eq-con-e}
\end{equation}%
and
\begin{equation}
\widetilde{E}\left( X,Y\right) =g\left( \widetilde{e}\,X,Y\right)
,\qquad X,Y\in {\mathfrak X}(M).  \label{eq-con-E}
\end{equation}%
Then
\begin{equation}
\widetilde{E}=\widetilde{S}-\frac{\widetilde{r}}{n}\,g.  \label{eq-con-barE}
\end{equation}%
The $\left( 0,2\right) $ tensor $\widetilde{E}$ is called tensor of Einstein
\cite{Stavre-82} with respect to the metric connection $\widetilde{\nabla }$%
. If $\widetilde{S}$ is symmetric then $\widetilde{E}$ is also symmetric.

\begin{defn-new}
A metric connection $\widetilde{\nabla }$ with symmetric Ricci tensor $%
\widetilde{S}$ will be called a {\bf Ricci-symmetric metric connection}.
\end{defn-new}

\begin{ex-new}
In a semi-Riemannian manifold $\left( M,g\right) $, a semi-symmetric metric
connection $\overline{\nabla }$ of K. Yano \cite{Yano-70}\ is given by
\[
\overline{\nabla }_{X}Y=\nabla _{X}Y+u(Y)X-g(X,Y)U,\qquad X,Y\in {\mathfrak X}%
(M),
\]%
where $\nabla $ is Levi-Civita connection, $U$ is a vector field and $u$ is
its associated $1$-form given by $u(X)=g(X,U)$. The Ricci tensor $\overline{S%
}$ with respect to $\overline{\nabla }$ is given by
\[
\overline{S}=S-(n-2)\alpha -{\rm trace}(\alpha )\,g,
\]%
where $S$ is the Ricci tensor and $\alpha $ is a $(0,2)$ tensor field
defined by%
\[
\alpha (X,Y)=(\nabla _{X}u)(Y)-u(X)u(Y)+\frac{1}{2}u(U)g(X,Y),\qquad
X,Y\in {\mathfrak X}(M).
\]%
The Ricci tensor $\overline{S}$ is symmetric if $1$-form $u$ is closed.
\end{ex-new}

\begin{ex-new}
An $(\varepsilon )$-almost para contact metric manifold $(M,\varphi ,\xi
,\eta ,g,\varepsilon )$ is given by
\[
\varphi ^{2}=I-\eta \otimes \xi ,\;\eta (\xi )=1,\;g\left( \varphi X,\varphi
Y\right) =g\left( X,Y\right) -\varepsilon \eta (X)\eta \left( Y\right) ,
\]%
where $\varphi $ is a tensor field of type $(1,1)$, $\eta $ is $1$-form, $%
\xi $ is a vector field and $\varepsilon =\pm 1$. An $(\varepsilon )$-almost
para contact metric manifold satisfying%
\[
\left( \nabla _{X}\varphi \right) Y=-\,g(\varphi X,\varphi Y)\xi
-\varepsilon \eta \left( Y\right) \varphi ^{2}X
\]%
is called an $(\varepsilon )$-para Sasakian manifold \cite{TKYK-09}. In an $%
(\varepsilon )$-para Sasakian manifold, the semi-symmetric metric connection
$\overline{\nabla }$ given by
\[
\overline{\nabla }_{X}Y=\nabla _{X}Y+\eta (Y)X-g(X,Y)\xi
\]%
is a Ricci symmetric metric connection.
\end{ex-new}

\begin{ex-new}
An almost contact metric manifold $(M,\varphi ,\xi ,\eta ,g)$ is given by%
\[
\varphi ^{2}=-I+\eta \otimes \xi ,\;\eta (\xi )=1,\;g\left( \varphi
X,\varphi Y\right) =g\left( X,Y\right) -\eta (X)\eta \left( Y\right) ,
\]%
where $\varphi $ is a tensor field of type $(1,1)$, $\eta $ is $1$-form and $%
\xi $ is a vector field. An almost contact metric manifold is a {\em %
Kenmotsu manifold} \cite{Kenmotsu-72} if
\[
\left( \nabla _{X}\varphi \right) Y=\,g(\varphi X,Y)\xi -\eta \left(
Y\right) \varphi X,
\]%
and is a Sasakian manifold \cite{Sasaki-60} if
\[
\left( \nabla _{X}\varphi \right) Y=g(X,Y)\xi -\eta \left( Y\right) X.
\]%
In an almost contact metric manifold $M$ the semi-symmetric metric
connection $\overline{\nabla }$ given by
\[
\overline{\nabla }_{X}Y=\nabla _{X}Y+\eta (Y)X-g(X,Y)\xi
\]%
is a Ricci symmetric metric connection if $M$ is Kenmotsu, but the
connection fails to be Ricci symmetric if $M$ is Sasakian.
\end{ex-new}

Let $(M,g)$ be an $n$-dimensional semi-Riemannian manifold equipped with a
metric connection $\widetilde{\nabla }$. A symmetric $\left( 0,2\right) $
tensor field $H$, which is covariantly constant with respect to $\widetilde{%
\nabla }$, is called a {\em special quadratic first integral} (for brevity
SQFI) \cite{Levine-Katzin-68-Tensor-1} with respect to $\widetilde{\nabla }$%
. The semi-Riemannian metric $g$ is always an SQFI. A set of SQFI tensors $%
\left\{ H_{1},\ldots ,H_{\ell }\right\} $ with respect to $\widetilde{\nabla
}$ is said to be {\em linearly independent} if
\[
c_{1}H_{1}+\cdots +c_{\ell }H_{\ell }=0,\qquad c_{1},\ldots ,c_{\ell }\in
{\bf R,}
\]%
implies that
\[
c_{1}=\cdots =c_{\ell }=0.
\]%
The set $\left\{ H_{1},\ldots ,H_{\ell }\right\} $ is said to be a complete
set if any SQFI tensor $H$ with respect to $\widetilde{\nabla }$ can be
written as
\[
H=c_{1}H_{1}+\cdots +c_{\ell }H_{\ell }\,,\qquad c_{1},\ldots ,c_{\ell }\in
{\bf R.}
\]%
The {\bf index }\cite{Thomas-39} of the manifold $M$ with respect to $%
\widetilde{\nabla }$, denoted by $i_{\widetilde{\nabla }}$, is
defined to be the number $\ell $ of members in a complete set
$\left\{ H_{1},\ldots ,H_{\ell }\right\} $.

We shall need the following Lemma:

\begin{lem}
\label{lem-con-1} Let $\left( M,g\right) $ be an $n$-dimensional
semi-Riemannian manifold equipped with a Ricci symmetric metric connection $%
\widetilde{\nabla }$. Then the following statements are true:

\begin{enumerate}
\item[{\rm (a)}] If $\widetilde{\nabla }_{X}\widetilde{S}=0$, then $%
\widetilde{\nabla }_{X}\widetilde{E}=0$. Conversely, if $\,\widetilde{r}$ is
constant and $\widetilde{\nabla }_{X}\widetilde{E}=0$ then $\widetilde{%
\nabla }_{X}\widetilde{S}=0$.

\item[{\rm (b)}] If $\widetilde{\nabla }_{X}\widetilde{S}\neq 0$ and $\psi $
is a non-vanishing differentiable function such that $\psi \widetilde{\nabla
}_{X}\widetilde{S}$ and $g$ are linearly dependent, then $\widetilde{\nabla }%
_{X}\widetilde{E}=0$.
\end{enumerate}
\end{lem}

The proof is similar to Lemmas 1.2 and 1.3 in \cite{Stavre-Smaranda-82} for
a semi-symmetric metric connection, and is therefore omitted.

\section{Quasi-conformal curvature tensor\label{sec-qcc-flat}}

Let $(M,g)$ be an $n$-dimensional $(n>3)$ semi-Riemannian manifold equipped
with a metric connection $\widetilde{\nabla }$. The conformal curvature
tensor $\widetilde{{\cal C}}$ with respect to the $\widetilde{\nabla }$ is
defined by \cite[p.~90]{Eisenhart-49}
\begin{eqnarray}
\widetilde{{\cal C}}(X,Y,Z,V) &=&\widetilde{R}\left( X,Y,Z,V\right) -\frac{1%
}{n-2}\left( \widetilde{S}\left( Y,Z\right) g(X,V)-\widetilde{S}\left(
X,Z\right) g(Y,V)\right.  \nonumber \\
&&+\ \left. g\left( Y,Z\right) \widetilde{S}\left( X,V\right) -g\left(
X,Z\right) \widetilde{S}\left( Y,V\right) \right)  \nonumber \\
&&+\ \frac{\widetilde{r}}{(n-1)(n-2)}\left( g\left( Y,Z\right) g\left(
X,V\right) -g\left( X,Z\right) g\left( Y,V\right) \right) ,
\label{eq-quasi-1a}
\end{eqnarray}%
and the concircular curvature tensor $\widetilde{{\cal Z}}$ with respect to $%
\widetilde{\nabla }$ is defined by (\cite{Yano-40}, \cite[p. 87]%
{Yano-Bochner-53})
\begin{equation}
\widetilde{{\cal Z}}\left( X,Y,Z,V\right) =\widetilde{R}\left(
X,Y,Z,V\right) -\frac{\widetilde{r}}{n\left( n-1\right) }\left( g\left(
Y,Z\right) g\left( X,V\right) -g\left( X,Z\right) g\left( Y,V\right) \right)
.  \label{eq-quasi-1b}
\end{equation}%
As a generalization of the notion of conformal curvature tensor and
concircular curvature tensor, the quasi-conformal curvature tensor $%
\widetilde{{\cal C}}_{\ast }$\ with respect to $\widetilde{\nabla }$ is
defined by \cite{Yano-Sawaki-68}
\begin{eqnarray}
\widetilde{{\cal C}}_{\ast }\left( X,Y,Z,V\right) &&=\, a\widetilde{R}\left(
X,Y,Z,V\right) +b\left( \widetilde{S}\left( Y,Z\right) g\left( X,V\right) -%
\widetilde{S}\left( X,Z\right) g\left( Y,V\right) \right.  \nonumber \\
&&+\ \left. g\left( Y,Z\right) \widetilde{S}\left( X,V\right) -g\left(
X,Z\right) \widetilde{S}\left( Y,V\right) \right)  \nonumber \\
&&-\ \frac{\widetilde{r}}{n}\left\{ \frac{a}{n-1}+2b\right\} \left( g\left(
Y,Z\right) g\left( X,V\right) -g\left( X,Z\right) g\left( Y,V\right) \right)
,  \label{eq-quasi-1}
\end{eqnarray}%
where $a$ and $b$ are constants. In fact, we have
\begin{equation}
\widetilde{{\cal C}}_{\ast }\left( X,Y,Z,V\right) =-(n-2)b\,\widetilde{{\cal %
C}}\left( X,Y,Z,V\right) +(a+(n-2)b)\widetilde{{\cal Z}}\left(
X,Y,Z,V\right) .  \label{eq-quasi-1c}
\end{equation}%
Since, there is no restrictions for manifolds if $a=0$ and $b=0$, therefore
it is essential for us to consider the case of $a\neq 0$ or $b\neq 0$. From (%
\ref{eq-quasi-1c}) it is clear that if $a=1$ and $b=-\,1/\left( n-2\right) $%
, then $\widetilde{{\cal C}}_{\ast }=\widetilde{{\cal C}}$; and if
$a=1$ and $b=0$, then $\widetilde{{\cal C}}_{\ast }=\widetilde{{\cal
Z}}$. \medskip

Now, we need the following:

\begin{defn-new}
A semi-Riemannian manifold $(M,g)$ equipped with a metric connection $%
\widetilde{\nabla }$ is said to be

\begin{enumerate}
\item[{\rm (a)}] $\widetilde{\nabla }$-quasi-conformally flat if $\widetilde{%
{\cal C}}_{\ast }=0$.

\item[{\rm (b)}] $\widetilde{\nabla }$-conformally flat if $\widetilde{{\cal %
C}}=0$.

\item[{\rm (c)}] $\widetilde{\nabla }$-concircularly flat if $\widetilde{%
{\cal Z}}=0$.
\end{enumerate}

\noindent In particular, with respect to the Levi-Civita connection $\nabla $%
, $\widetilde{\nabla }$-quasi-conformally flat, $\widetilde{\nabla }$%
-conformally flat and $\widetilde{\nabla }$-concircularly flat become simply
quasi-conformally flat, conformally flat and concircularly flat respectively.
\end{defn-new}

\begin{defn-new}
A semi-Riemannian manifold $(M,g)$ equipped with a metric connection $%
\widetilde{\nabla }$ is said to be

\begin{enumerate}
\item[{\rm (a)}] $\widetilde{\nabla }$-quasi-conformally symmetric if $%
\widetilde{\nabla }\,\widetilde{{\cal C}}_{\ast }=0$.

\item[{\rm (b)}] $\widetilde{\nabla }$-conformally symmetric if $\widetilde{%
\nabla }\widetilde{{\cal C}}=0$.

\item[{\rm (c)}] $\widetilde{\nabla }$-concircularly symmetric if $%
\widetilde{\nabla }\widetilde{{\cal Z}}=0$.
\end{enumerate}

\noindent In particular, with respect to the Levi-Civita connection $\nabla $%
, $\widetilde{\nabla }$-quasi-conformally symmetric, $\widetilde{\nabla }$%
-conformally symmetric and $\widetilde{\nabla }$-concircularly symmetric
become simply quasi-conformally symmetric, conformally symmetric and
concircularly symmetric respectively.
\end{defn-new}

\begin{th}
\label{th-quasi} Let $M$ be a semi-Riemannian manifold of dimension $n>2$. Then $M$ is $\widetilde{\nabla }$-quasi-conformally flat
if and only if one of the following statements is true{\rm :}

\begin{description}
\item[{\rm (i)}] $a+(n-2)b=0$, $a\neq 0\neq b$ and $M$ is $\widetilde{\nabla
}$-conformally flat.

\item[{\rm (ii)}] $a+(n-2)b\not=0$, $a\neq 0$, $M$ is $\widetilde{\nabla }$%
-conformally flat and $\widetilde{\nabla }$-concircularly flat.

\item[{\rm (iii)}] $a+(n-2)b\not=0$, $a=0$ and Ricci tensor $\widetilde{S}$
with respect to $\widetilde{\nabla }$ satisfies%
\begin{equation}
\widetilde{S}-\frac{\widetilde{r}}{n}\,g=0,  \label{eq-S}
\end{equation}%
where $\widetilde{r}$ is the scalar curvature with respect to $\widetilde{%
\nabla}$.
\end{description}
\end{th}
\noindent {\bf Proof.}
Using $\widetilde{{\cal C}}_{\ast }=0$ in (\ref%
{eq-quasi-1}) we get
\begin{eqnarray}
0 &=&a\widetilde{R}\left( X,Y,Z,V\right) +b\left( \widetilde{S}\left(
Y,Z\right) g\left( X,V\right) -\widetilde{S}\left( X,Z\right) g\left(
Y,V\right) \right.  \nonumber \\
&&+\ \left. g\left( Y,Z\right) \widetilde{S}\left( X,V\right) -g\left(
X,Z\right) \widetilde{S}\left( Y,V\right) \right)  \nonumber \\
&&-\ \frac{\widetilde{r}}{n}\left( \frac{a}{n-1}+2b\right) \left( g\left(
Y,Z\right) g\left( X,V\right) -g\left( X,Z\right) g\left( Y,V\right) \right)
,  \label{eq-quasi-2}
\end{eqnarray}%
from which we obtain
\begin{equation}
(a+(n-2)b)\left( \widetilde{S}-\frac{\widetilde{r}}{n}\,g\right) =0.
\label{eq-quasi-2a}
\end{equation}

\noindent {\bf Case 1.} $a+(n-2)b=0$ and $a\neq 0\neq b$. Then from (\ref{eq-quasi-1})
and (\ref{eq-quasi-1a}), it follows that $(n-2)b\,\widetilde{{\cal C}}=0$,
which gives $\widetilde{{\cal C}}=0$. This gives the statement (i).

\noindent{\bf Case 2.} $a+(n-2)b\not=0$ and $a\neq 0$. Then from (\ref{eq-quasi-2a})
\begin{equation}
\widetilde{S}\left( Y,Z\right) =\frac{\widetilde{r}}{n}g(Y,Z).
\label{eq-quasi-s}
\end{equation}%
Using (\ref{eq-quasi-s}) in (\ref{eq-quasi-2}), we get
\begin{equation}
a(\widetilde{R}\left( X,Y,Z,V\right) -\frac{\widetilde{r}}{n(n-1)}\left(
g(Y,Z)g(X,V)-g(X,Z)g(Y,V)\right) )=0.  \label{eq-quasi-2b}
\end{equation}%
Since $a\not=0$, then by (\ref{eq-quasi-1b}) $\widetilde{{\cal Z}}=0$ and by
using (\ref{eq-quasi-2b}), (\ref{eq-quasi-s}) in (\ref{eq-quasi-1a}), we get
$\widetilde{{\cal C}}=0$. This gives the statement (ii).

\noindent{\bf Case 3.} $a+(n-2)b\not=0$ and $a=0$, we get (\ref{eq-S}). This gives
the statement (iii). Converse is true in all cases.

\begin{cor}
\label{cor-qcf} {\rm \cite[Theorem 5.1]{Tripathi-Gupta-10}} Let $M$ be a semi-Riemannian manifold of dimension $n>2$. Then $M$ is quasi-conformally flat if and only if one of
the following statements is true{\rm :}

\begin{description}
\item[{\rm (i)}] $a+(n-2)b=0$, $a\neq 0\neq b$ and $M$ is conformally flat.

\item[{\rm (ii)}] $a+(n-2)b\not=0$, $a\neq 0$, $M$ is of constant curvature.

\item[{\rm (iii)}] $a+(n-2)b\not=0$, $a=0$ and $M$ is Einstein manifold.
\end{description}
\end{cor}

\begin{rem-new}
In \cite{Amur-Mar-77-Tensor}, the following three results are known:

\begin{enumerate}
\item[{\rm (a)}] \cite[Proposition 1.1]{Amur-Mar-77-Tensor} A
quasi-conformally flat manifold is either conformally flat or Einstein.

\item[{\rm (b)}] \cite[Corollary 1.1]{Amur-Mar-77-Tensor} A
quasi-conformally flat manifold is conformally flat if the constant $a\not=0$%
.

\item[{\rm (c)}] \cite[Corollary 1.2]{Amur-Mar-77-Tensor} A
quasi-conformally flat manifold is Einstein if the constants $a=0$ and $%
b\not=0$.
\end{enumerate}

However, the converses need not be true in these three results. But, in
Corollary~\ref{cor-qcf} we get a complete classification of
quasi-conformally flat manifolds.
\end{rem-new}

\section{$\widetilde{\protect\nabla }$-Quasi-conformally symmetric manifolds
\label{sec-qcc-index}}

Let $(M,g)$ be an $n$-dimensional semi-Riemannian manifold equipped with the
metric connection $\widetilde{\nabla }$. Let $\widetilde{R}$ be the
curvature tensor of $M$ with respect to the metric connection $\widetilde{%
\nabla }$. If $H$ is a parallel symmetric $(0,2)$ tensor with respect to the
metric connection $\widetilde{\nabla }$, then we easily obtain
\begin{equation}
H((\widetilde{\nabla }_{U}\widetilde{R})(X,Y)Z,V)+H(Z,(\widetilde{\nabla }%
_{U}\widetilde{R})(X,Y)V)=0,\qquad X,Y,Z,V,U\in {\mathfrak X}(M).
\label{eq-con-tens}
\end{equation}%
The solutions $H$ of (\ref{eq-con-tens}) is closely related to the index of
quasi-conformally symmetric and concircularly symmetric manifold with
respect to the $\widetilde{\nabla }$.

\begin{lem}
\label{lem-quasi-1} Let $\left( M,g\right) $ be an $n$-dimensional
semi-Riemannian $\widetilde{\nabla }$-quasi-conformally symmetric manifold, $%
n>2$ and $b\neq 0$. Then
\begin{equation}
{\rm trace}(\widetilde{\nabla }_{U}\widetilde{E})=0.  \label{eq-quasi-4}
\end{equation}
\end{lem}

\noindent {\bf Proof.}
Using (\ref{eq-con-barE}) in (\ref{eq-quasi-1}) we get
\begin{eqnarray}
\widetilde{{\cal C}}_{\ast }\left( X,Y,Z,V\right) &=&a\widetilde{R}\left(
X,Y,Z,V\right) +b\left( \widetilde{E}\left( Y,Z\right) g\left( X,V\right) -%
\widetilde{E}\left( X,Z\right) g\left( Y,V\right) \right.  \nonumber \\
&&+\ \left. g\left( Y,Z\right) \widetilde{E}\left( X,V\right) -g\left(
X,Z\right) \widetilde{E}\left( Y,V\right) \right)  \nonumber \\
&&-\ \frac{a\ \widetilde{r}}{n\left( n-1\right) }\left( g\left( Y,Z\right)
g\left( X,V\right) -g\left( X,Z\right) g\left( Y,V\right) \right) .
\label{eq-quasi-3}
\end{eqnarray}%
Taking covariant derivative of (\ref{eq-quasi-3}) and using $\widetilde{%
\nabla }_{U}\widetilde{{\cal C}}_{\ast }=0$, we get
\begin{eqnarray}
a(\widetilde{\nabla }_{U}\widetilde{R})\left( X,Y,Z,V\right) &&=b\left( (%
\widetilde{\nabla }_{U}\widetilde{E})\left( X,Z\right) g\left( Y,V\right) -(%
\widetilde{\nabla }_{U}\widetilde{E})\left( Y,Z\right) g\left( X,V\right)
\right.  \nonumber \\
&&-\ \left. g\left( Y,Z\right) (\widetilde{\nabla }_{U}\widetilde{E})\left(
X,V\right) +g\left( X,Z\right) (\widetilde{\nabla }_{U}\widetilde{E})\left(
Y,V\right) \right)  \nonumber \\
&&+\ \frac{a(\widetilde{\nabla }_{U}\widetilde{r})}{n\left( n-1\right) }%
\left( g\left( Y,Z\right) g\left( X,V\right) - g\left( X,Z\right) g\left(
Y,V\right) \right) .  \label{eq-quasi-6}
\end{eqnarray}%
Contracting (\ref{eq-quasi-6}) with respect to $Y${\rm \ }and $Z$ and using (%
\ref{eq-curv-con-semi-hay-1}) and (\ref{eq-curv-con-semi-hay-2}), we get
\begin{eqnarray}
a(\widetilde{\nabla }_{U}\widetilde{S})\left( X,V\right)
&=&-\,b\,{\rm trace}(\widetilde{\nabla }_{U}\widetilde{E}) g\left( X,V\right)  \nonumber \\
&&-\,\left( n-2\right) b(\widetilde{\nabla }_{U}\widetilde{E})\left(
X,V\right) +\frac{a(\widetilde{\nabla }_{U}\widetilde{r})}{n}g\left(
X,V\right) .  \label{eq-quasi-9}
\end{eqnarray}
Using (\ref{eq-quasi-9}), we
get (\ref{eq-quasi-4}).

\begin{th}
\label{th-quasi-1} If $\left( M,g\right) $ is an $n$-dimensional
semi-Riemannian $\widetilde{\nabla }$-quasi-conformally symmetric manifold, $%
n>2$ and $b\neq 0$, then the equation $(\ref{eq-con-tens})$ takes the form
\begin{equation}
\det \left(
\begin{array}{cc}
H\left( X,Z\right) -\displaystyle\frac{1}{n}\,{\rm trace}\left( H\right)
g\left( X,Z\right) & \quad H\left( Y,V\right) -\displaystyle\frac{1}{n}\,%
{\rm trace}\left( H\right) g\left( Y,V\right) \medskip \\
(\widetilde{\nabla }_{U}\widetilde{E})\left( X,Z\right) & (\widetilde{\nabla
}_{U}\widetilde{E})\left( Y,V\right)%
\end{array}%
\right) =0.  \label{eq-quasi-14}
\end{equation}%
If $\widetilde{\nabla }_{U}\widetilde{E}\neq 0$, then $(\ref{eq-quasi-14})$
has the general solution
\begin{equation}
H_{U}\left( X,Y\right) =f(\widetilde{\nabla }_{U}\widetilde{S})\left(
X,Y\right) +\frac{1}{n}\left( {\rm trace}(H_{U})-f(\widetilde{\nabla }_{U}%
\widetilde{r})\right) g\left( X,Y\right) ,  \label{eq-quasi-tens-1}
\end{equation}%
where $f$ is an arbitrary non-vanishing differentiable function.
\end{th}

\noindent {\bf Proof.}
Using (\ref{eq-quasi-6}) in (\ref{eq-con-tens}), we get
\begin{eqnarray}
0 &=&b\left( (\widetilde{\nabla }_{U}\widetilde{E})\left( X,Z\right) H\left(
Y,V\right) -(\widetilde{\nabla }_{U}\widetilde{E})\left( Y,Z\right) H\left(
X,V\right) \right.  \nonumber \\
&&-\ g\left( Y,Z\right) H((\widetilde{\nabla }_{U}\widetilde{e})X,V)+g\left(
X,Z\right) H((\widetilde{\nabla }_{U}\widetilde{e})Y,V)  \nonumber \\
&&+\ (\widetilde{\nabla }_{U}\widetilde{E})\left( X,V\right) H\left(
Y,Z\right) -(\widetilde{\nabla }_{U}\widetilde{E})\left( Y,V\right) H\left(
X,Z\right)  \nonumber \\
&&-\ \left. g\left( Y,V\right) H((\widetilde{\nabla }_{U}\widetilde{e}%
)X,Z)+g\left( X,V\right) H((\widetilde{\nabla }_{U}\widetilde{e})Y,Z)\right)
\nonumber \\
&&+\ \frac{a(\widetilde{\nabla }_{U}\widetilde{r})}{n\left( n-1\right) }%
\left( g\left( Y,Z\right) H\left( X,V\right) -g\left( X,Z\right) H\left(
Y,V\right) \right.  \nonumber \\
&&+\ \left. g\left( Y,V\right) H\left( X,Z\right) -g\left( X,V\right)
H\left( Y,Z\right) \right) .  \label{eq-quasi-tens-3}
\end{eqnarray}%
Let $\{e_{1},\ldots ,e_{n}\}$ be an orthonormal basis
of vector fields in $M$. Taking $X=Z=e_{i}$ in (\ref{eq-quasi-tens-3}) and summing upto $n$ terms, then using (\ref{eq-quasi-4}), we have
\begin{eqnarray}
0 &=&b\left( \left( n-1\right) H((\widetilde{\nabla }_{U}\widetilde{e}%
)Y,V)+H((\widetilde{\nabla }_{U}\widetilde{e})V,Y)\right.  \nonumber \\
&&-\left. {\rm trace}(H)(\widetilde{\nabla }_{U}\widetilde{E})\left(
Y,V\right) -g\left( Y,V\right) \sum_{i=1}^{n}H((\widetilde{\nabla }_{U}%
\widetilde{e})e_{i},e_{i})\right)  \nonumber \\
&&+\ \frac{a(\widetilde{\nabla }_{U}\widetilde{r})}{n\left( n-1\right) }%
\left( {\rm trace}(H)g\left( Y,V\right) -nH\left( Y,V\right) \right) .
\label{eq-quasi-tens-6}
\end{eqnarray}%
Interchanging $Y$ and $V$ in (\ref{eq-quasi-tens-6}) and subtracting the so
obtained formula from (\ref{eq-quasi-tens-6}), we deduce that
\begin{equation}
H((\widetilde{\nabla }_{U}\widetilde{e})Y,V)=H((\widetilde{\nabla }_{U}%
\widetilde{e})V,Y).  \label{eq-quasi-tens-9}
\end{equation}%
Now, interchanging $X$ and $Z$, $Y$ and $V$ in (\ref{eq-quasi-tens-3}) and
taking the sum of the resulting equation and (\ref{eq-quasi-tens-3}) and
using (\ref{eq-quasi-tens-6}) and (\ref{eq-quasi-tens-9}), we get (\ref%
{eq-quasi-14}). If $\widetilde{\nabla }_{U}\widetilde{E}\neq 0$, then using (%
\ref{eq-con-barE}) leads to (\ref{eq-quasi-tens-1}).

\begin{th}
\label{th-conformal-tens-1} If $\left( M,g\right) $ is an $n$-dimensional
semi-Riemannian $\widetilde{\nabla }$-quasi-conformally symmetric manifold, $%
n>2$ and $b\neq 0$, and if there is a vector field $U$ so that
\begin{equation}
\widetilde{\nabla }_{U}\widetilde{E}=0\quad {\rm and}\quad \widetilde{\nabla
}_{U}\widetilde{r}\neq 0,  \label{eq-quasi-tens-22}
\end{equation}%
then the solution of a equation $(\ref{eq-con-tens})$ is $H=f\,g$, where $f$
is a differentiable non-vanishing function.
\end{th}

\noindent {\bf Proof.}
Using (\ref{eq-quasi-tens-22}), (\ref{eq-quasi-tens-3}%
) becomes
\begin{equation}
g\left( Y,Z\right) H\left( X,V\right) -g\left( X,Z\right) H\left( Y,V\right)
+g\left( Y,V\right) H\left( X,Z\right) -g\left( X,V\right) H\left(
Y,Z\right) =0,  \label{eq-quasi-tens-23}
\end{equation}%
Interchanging $X$ and $Z$, $Y$ and $V$ in (\ref{eq-quasi-tens-23}) and
taking the sum of the resulting equation and (\ref{eq-quasi-tens-23}), we
get
\[
g\left( X,Z\right) H\left( Y,V\right) -g\left( Y,V\right) H\left( X,Z\right)
=0.
\]%
Therefore the tensor fields $H$ and $g$ are proportional.

\begin{th}
\label{th-quasi-index-1} Let $\left( M,g\right) $ be an $n$-dimensional
semi-Riemannian $\widetilde{\nabla }$-quasi-conformally symmetric manifold, $%
n>2$ and $b\neq 0$. If there is a vector field $U$ satisfying the condition $%
(\ref{eq-quasi-tens-22})$, then $i_{\widetilde{\nabla }}=1$.
\end{th}

\noindent {\bf Proof.}
By Theorem \ref{th-conformal-tens-1} and from the fact that
$\widetilde{\nabla }_{U}g=0$ and $\widetilde{\nabla }_{U}H=0$, it
follows that $f$ is constant. Thus $i_{\widetilde{\nabla }}=1$.

\begin{th}
\label{th-quasi-index-2} Let $\left( M,g\right) $ be an $n$-dimensional
semi-Riemannian $\widetilde{\nabla }$-quasi-conformally symmetric manifold, $%
n>2$ and $b\neq 0$, for which the tensor field $\widetilde{E}$ is not
covariantly constant with respect to the Ricci symmetric metric connection $%
\widetilde{\nabla }$. If $i_{\widetilde{\nabla }}>1$, then there is a vector
field $U$, so that the equation
\begin{equation}
\widetilde{\nabla }_{U}H=0  \label{eq-quasi-index-1}
\end{equation}%
has the fundamental solutions
\begin{equation}
H_{1}=g,\qquad H_{2}=\psi \widetilde{\nabla }_{U}\widetilde{S},
\label{eq-quasi-index-2}
\end{equation}%
where $\psi $ is a differentiable non-vanishing function.
\end{th}

\noindent {\bf Proof.}
Given that $\widetilde{\nabla }_{U}E\neq 0$, there is $U$ so
that the tensorial equation (\ref{eq-con-tens}) has general solution
which depends on $U$. $g$ is obviously a solution of
(\ref{eq-quasi-index-1} ) because $\widetilde{\nabla }_{U}g=0$, $g$
is also satisfy the tensorial equation (\ref{eq-con-tens}) and
$H_{U}$ given by (\ref{eq-quasi-tens-1}) is also a solution of
(\ref{eq-quasi-index-1}). Equation (\ref{eq-quasi-index-1} ) has at
least two solution as $i_{\widetilde{\nabla }}>1$. These two
solution are independent.By Lemma~\ref{lem-con-1}(b) $\psi
\widetilde{\nabla }_{U}\widetilde{S}$ and $g$ are independent and we
get two fundamental
solution of $\widetilde{\nabla }_{U}\widetilde{H}=0$ which is $%
H_{1}=g,\,H_{2}=\psi \widetilde{\nabla }_{U}\widetilde{S}$, where
$\psi $ is a differentiable non-vanishing function.

\begin{th}
\label{th-quasi-index-3} Let $\left( M,g\right) $ be an $n$-dimensional
semi-Riemannian $\widetilde{\nabla }$-quasi-conformally symmetric manifold, $%
n>2$ and $b\neq 0$, for which the tensor field $\widetilde{E}$ is not
covariantly constant with respect to the metric connection $\widetilde{%
\nabla }$. Then $1\leq i_{\widetilde{\nabla }}\leq n+1$.
\end{th}
\noindent {\bf Proof.}
Let $U_{i}$, $i=1,\ldots ,p$ be independent vector fields, for which
\[
\widetilde{\nabla }_{U_{i}}\widetilde{E}\neq 0,
\]%
and let $\psi _{i}{\widetilde{\nabla }_{U_{i}}}\widetilde{S}$ and $g$ be the
fundamental solutions of ${\widetilde{\nabla }_{U_{i}}}\widetilde{H}=0$.
Obviously $p<n$, as $U_{i}$ are independent. Therefore we have $p+1$
solutions.
This completes the proof.

\begin{rem-new}
The previous results of this section will be true for $\widetilde{\nabla }$%
-conformally symmetric semi-Riemannian manifold, where $\widetilde{\nabla }$
is any Ricci symmetric metric connection.
\end{rem-new}

\begin{th}
\label{th-conc-1} If $\left( M,g\right) $ be an $n$-dimensional
semi-Riemannian $\widetilde{\nabla }$-concircularly symmetric manifold, then
the equation $(\ref{eq-con-tens})$ takes the form
\begin{equation}
\det \left(
\begin{array}{cc}
H\left( X,Z\right) & H\left( Y,V\right) \medskip \\
g\left( X,Z\right) & g\left( Y,V\right)%
\end{array}%
\right) =0.  \label{eq-conc-h}
\end{equation}
\end{th}
\noindent {\bf Proof.}
Taking covariant derivative of (\ref{eq-quasi-1b})
and using $\widetilde{\nabla }_{U}\widetilde{{\cal Z}}=0$, we get%
\[
(\widetilde{\nabla }_{U}\widetilde{R})\left( X,Y,Z,V\right) =\frac{%
\widetilde{\nabla }_{U}\widetilde{r}}{n\left( n-1\right) }\left( g\left(
Y,Z\right) g\left( X,V\right) -g\left( X,Z\right) g\left( Y,V\right) \right)
,
\]%
which, when used in (\ref{eq-con-tens}), yields%
\begin{eqnarray}
0 &=& \frac{\widetilde{\nabla }_{U}\widetilde{r}}{n\left( n-1\right) }\left(
g\left( Y,Z\right) H\left( X,V\right) -g\left( X,Z\right) H\left( Y,V\right)\right.
\nonumber\\&&\left.+\ g\left( Y,V\right) H\left( X,Z\right) - g\left( X,V\right) H\left(
Y,Z\right) \right) .  \label{eq-conc-tens-2}
\end{eqnarray}%
Now, we interchange $X$ with $Z$, and $Y$ with $V$ in
(\ref{eq-conc-tens-2}) and take the sum of the resulting equation
and (\ref{eq-conc-tens-2}), we get (\ref{eq-conc-h}).

\begin{th}
\label{th-conc-index-1} Let $\left(M,g\right)$ be an $n$-dimensional
semi-Riemannian $\widetilde{\nabla}$-concircularly symmetric manifold. Then $%
i_{\widetilde{\nabla}}=1$.
\end{th}

\noindent {\bf Proof.}
By Theorem \ref{th-conc-1} and from the fact that $%
\widetilde{\nabla }_{U}g=0$ and $\widetilde{\nabla }_{U}H=0$, we get
$i_{ \widetilde{\nabla }}=1$.

\section{Discussion}

A semi-Riemannian manifold is said to be {\em decomposable} \cite{Thomas-39}
(or locally reducible) if there always exists a local coordinate system $%
\left( x^{i}\right) $ so that its metric takes the form
\[
ds^{2}=\sum_{a,b=1}^{r}g_{ab}dx^{a}dx^{b}+\sum_{\alpha ,\beta
=r+1}^{n}g_{\alpha \beta }dx^{\alpha }dx^{\beta },
\]
where $g_{ab}$ are functions of $x^{1},\ldots ,x^{r}$ and $g_{\alpha \beta }$
are functions of $x^{r+1},\ldots ,x^{n}$. A semi-Riemannian manifold is said
to be {\em reducible} if it is isometric to the product of two or more
semi-Riemannian manifolds; otherwise it is said to be {\em irreducible} \cite%
{Thomas-39}. A reducible semi-Riemannian manifold is always decomposable but
the converse need not be true.

The concept of the index of\ a (semi-)Riemannian manifold gives a striking
tool to decide the reducibility and decomposability of (semi-)Riemannian
manifolds. For example, a Riemannian manifold is decomposable if and only if
its index is greater than one \cite{Thomas-39}. Moreover, a complete
Riemannian manifold is reducible if and only if its index is greater than
one \cite{Thomas-39}. A second order $(0,2)$-symmetric parallel tensor is
also known as a special Killing tensor of order two. Thus, a Riemannian
manifold admits a special Killing tensor other than the Riemannian metric $g$
if and only if the manifold is reducible \cite{Eisenhart-23}, that is the
index of the manifold is greater than $1$. In 1951, E.M. Patterson \cite%
{Patterson-51} found a similar result for semi-Riemannian manifolds. In
fact, he proved that a semi-Riemannian manifold $\left( M,g\right) $
admitting a special Killing tensor $K_{ij}$, other than $g$, is reducible if
the matrix $(K_{ij})$ has at least two distinct characteristic roots at
every point of the manifold. In this case, the index of the manifold is
again greater than $1$.

By Theorem \ref{th-quasi-index-3}, we conclude that a $\widetilde{\nabla }$%
-quasi-conformally symmetric Riemannian manifold (where $\widetilde{\nabla }$
is any Ricci symmetric metric connection, not necessarily Levi-Civita
connection) is decomposable and it is reducible if the manifold is complete.

It is known that the maximum number of linearly independent Killing
tensors of order $2$ in a semi-Riemannian manifold $\left(
M^{n},g\right) $ is ${\frac{1}{12}n(n+1)^{2}(n+2)}$, which is
attained if and only if $M$ is of constant curvature. The maximum
number of linearly independent Killing tensors in a four dimensional
spacetime is $50$ and this number is attained if and only if the
spacetime is of constant curvature \cite{Hauser-Malhiot-75}. But,
from Theorem \ref{th-quasi-index-3}, we also conclude that the
maximum number of linearly independent special Killing tensors in a
$4$-dimensional Robertson-Walker spacetime \cite[p.~341]{ONeill-83}
is $5$.

\noindent Department of Mathematics, Faculty of Science,\\ Banaras
Hindu University, Varanasi 221005, India\\Email:{\em
mmtripathi66@yahoo.com}\medskip

\noindent Department of Mathematics, Faculty of Science,\\ Banaras
Hindu University, Varanasi 221005, India\\Email:{\em punam\_2101@yahoo.co.in}\medskip

\noindent {Department of Applied Mathematics,\\ Chonnam National
University, Yosu 550-749, Korea\\
Email:{\em jskim0807@yahoo.co.kr}
\end{document}